# Multiple 3-Coloring, an Approach to 4-Coloring of Planar Graphs

Li, Shaoqing

Anhui Jianzhu University, Hefei, Anhui, China.

**Abstract:** A planar graph can be embedded in a piecewise linear manifold, and the lattice on each linear piece can be colored with 3-coloring. If a planar graph can be colored with multiple 3-coloring, i.e. coloring the graph in pieces with different 3-color subsets of 4 colors, then the graph is 4-colorable. In this paper, multiple 3-coloring was introduced, and then the combination and partition of planar graphs for multiple 3-coloring was studied. The study reveals that planar graphs can generally be decomposed into independent subgraphs, and each subgraph can be triangulated into a symmetric structure for multiple 3-coloring.

## 1. Introduction

According to the four-color theorem, any planar map can be colored with four colors so that any adjacent regions have different colors. In graph theory, the coloring of maps is often described as the face coloring of polyhedral graphs, or the vertex coloring of maximal planar graphs.

First proposed in 1852 [1, 2], the four-color theorem has been focused on for more than one and a half hundred years. In 1879, Alfred Kempe presented a "proof" of the four-color theorem with the method of Kempe chains [3]. Although his proof was shown incorrect by Percy Heawood in 1890, his method was successfully used to prove the five-color theorem [1]. In 1976, Kenneth Appel and Wolfgang Haken announced that they had proved the four-color theorem with the aid of a computer [4], but many scientists are skeptical of the proof [2, 5, 6] as it hasn't figured out how to color a planar graph.

The key to the 4-coloring of planar graphs is the combination of odd elements. Elizabeth J. Hartung proved that the 12 pentagonal faces of a fullerene graph can be combined into 6 pairs (i.e. the nice pairs in [7]), and the fullerene can be mainly colored with face 3-coloring except for incompatible faces adjacent to the connecting chain of each pair [8].

A polyhedral graph can be embedded in a piecewise linear manifold with each non-hexagonal face corresponding to the vertex of a cone manifold. The lattice on the cone manifold to an odd face can be partitioned into three zones of linear pieces with Kempe chains, and each zone can be colored with a different 3-color subset of 4 colors. Such coloring is called multiple 3-coloring. In this paper, face coloring for the hexagonal lattice on plane, tube, and cone manifold was studied respectively and multiple 3-coloring was introduced. Then combination and partition for multiple 3-coloring of polyhedral graphs were studied.

The remainder of the paper is as follows. In section 2, regular colorings of the planar hexagonal lattice were studied, and congruence class and congruence chain in graph 3-coloring were introduced. In section 3, regular coloring for tubular hexagonal lattices was studied. In section 4, multiple 3-coloring together with its special case (3+1)-coloring was introduced to the hexagonal lattice around an odd face. In section 5, polyhedral graphs were decomposed into independent subgraphs for M3-coloring with combinations of odd faces. In section 6, the combination graph was introduced, then the partition and triangulation of subgraph for multiple 3-coloring was studied. In section 7, orbital coloring for symmetric graphs was studied. In section 8, the conclusion and outlook were presented.

## 2. Coloring a planar hexagonal lattice

### 2.1 Two regular colorings

There are two regular colorings for a hexagonal lattice. It is easier to illustrate with its dual, a triangular lattice. Color the vertices of an equilateral triangle with three colors, then, color the vertices mirrored by its edges. If each mirrored vertex has the same color as the original vertex, the result is a 3-coloring as Figure 1(a). If the mirrored vertex has a fourth color different from the original three, the result is a regular 4-coloring (R4-coloring for short) as Figure 1(b).

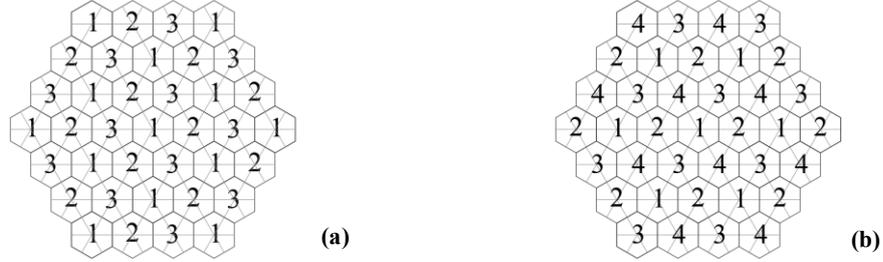

**Fig1. Two regular coloring of a hexagonal/triangular lattice.**
(a) 3-coloring.   (b) regular 4-coloring

The 3-coloring of the triangular lattice has symmetries of 120° rotation around each vertex. For a maximal planar graph, the cone angle around an even vertex is a multiple of 120°, so the triangular lattice around the vertex is 3-colorable. The cone angle around an odd vertex is not a multiple of 120°, and the triangular lattice around the vertex is not 3-colorable. Therefore, we can get the following theorem.

**Theorem 1.** (Saaty and Kainen [9] ). A maximal planar graph is vertex 3-colorable if and only if each vertex has an even degree.

The triangular lattice with R4-coloring has symmetries of 180° rotation around each vertex. The 180° cone angle corresponds to 3 degrees of the vertex. Then, we can get the following theorem.

**Theorem 2.** A maximal planar graph is vertex R4-colorable if and only if the degree of each vertex is a multiple of three.

Polyhedral graphs whose non-hexagonal faces are all triangles or squares are called (3,6)-fullerenes and (4,6)-fullerenes respectively. According to Theorem 1 and 2, the following corollaries can be obtained.

**Corollary 1.** Any (4,6)-fullerene is face 3-colorable.

**Corollary 2**. Any (3,6)-fullerene is face R4-colorable.

### 2.2 Congruence classes and chains in face 3-coloring

**Theorem 3.** For a hexagonal lattice, two faces have the same color in 3-coloring if and only if their relative Coxeter coordinates are congruent modulo 3.

**Proof.** In the 3-coloring, the face colors along a lattice line alternate in cycles of three colors, and two faces in a lattice line have the same color if and only if their distance is a multiple of 3. Just as in Figure 2(a), two faces connected with Coxeter coordinates ($p, q$) can be connected with another polyline path ($p, p+q$) which has a 120° turn at the joint. Because the lattice in 3-coloring has the symmetries of 120° rotation, the face colors alternate along the polyline the same as along a line. As the distance along the polyline is $2p + q = 3p - (p - q)$, so, if and only if the two Coxeter coordinates ($p, q$) are congruent modulo 3 the distance is a multiple of 3, and the two faces have the same color.

□

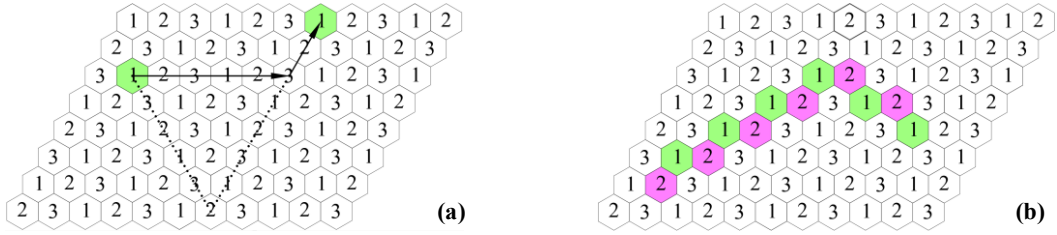

**Fig2. Congruence classes and congruence chains in face 3-coloring.**

According to the congruence of the Coxeter coordinates, the faces of a hexagonal lattice can be classified into three congruence classes, each having a different color in 3-coloring.

Successive adjacent faces in a congruence class form a congruence chain. If the chain does not cross itself, it is called a plain chain [8]. In this paper, all the congruence chains are plain chains by default. Just as in Figure 2(b), a congruence chain of color 1 and its adjacent congruence chain of color 2 form a Kempe chain. In this paper, all Kempe chains are such chains each composed of two congruence chains.

## 3. Coloring the hexagonal lattice on a tube

The hexagonal lattice on a tube can be unfolded along a linear congruence chain, just as in Figure 3. After unfolding, the Coxeter coordinates of two locations of a face on the chain, as *(p, q)* in Figure 3(a), are called the characteristic Coxeter coordinates of the tube.

**Theorem 4.** If two characteristic Coxeter coordinates of a tube are congruent modulo 3, the lattice is 3-colorable. If two parameters are not congruent modulo 3, the lattice can be colored mainly with three colors except for one linear congruence chain with a fourth color.

**Proof.** Just as in Figure 3(a), The hexagonal lattice on a tube can also be viewed as a series of consecutive linear chains. with 3-coloring, the faces on each chain have the same color, and adjacent chains have different colors. As each chain has a face on the dashed polyline, the number of chains can be counted along the polyline. As mentioned above, with 3-coloring the colors of faces along the polyline alternate in cycles of three colors. If two characteristic Coxeter coordinates of a tube are congruent modulo 3, the total number of linear chains is a multiple of 3, and the two locations of the unfolding chain have the same color, so the lattice is face 3-colorable.

Just as in Figure 3(b), if the total number of linear chains is of form $3n+1$ where *n* is a natural number, we can assign a fourth color to one linear chain, and then other faces can be colored with 3-coloring.

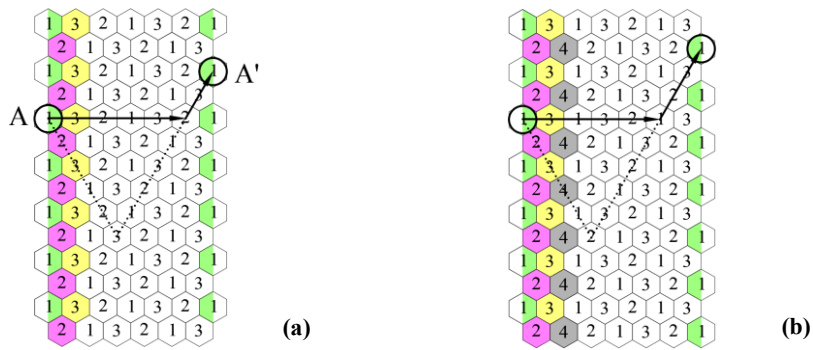

**Fig3. Face coloring for tubular hexagonal lattice. (a).** The lattice colored with 3-coloring. **(b).** The lattice colored with only one linear chain for a fourth color.

A tube with characteristic Coxeter coordinates ($p, q$) can be unfolded along linear chains in two different directions, just as in Figure 4. If unfolded as in Figure 4(a), the number of linear chains parallel to the unfolding chain is ($2p+q$). If unfolded as in Figure 4(b), the number is ($p+2q$). As ($2p+q$) + ($p+2q$) = $3(p+q)$, if $p$ and $q$ are not congruent modulo 3, one of the two numbers must be of form $3n+1$, and the corresponding unfolding can be colored as Figure 3(b).

□

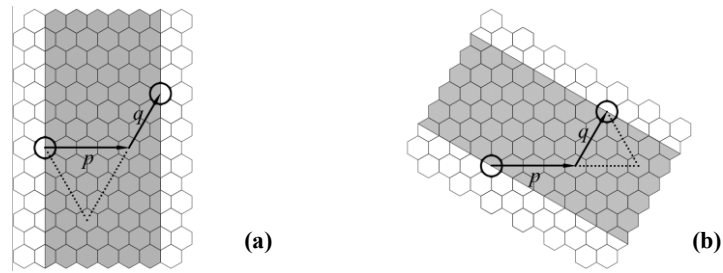

**Fig4. Unfolding a tubular hexagonal lattice through chains in different directions.**

## 4. Coloring around an odd face

The hexagonal lattice around a non-hexagonal face can be embedded in a cone manifold with the cone vertex corresponding to the non-hexagonal face. The neighborhood of an odd face is an odd cycle that needs at least three colors, so the lattice including an odd face needs at least 4 colors.

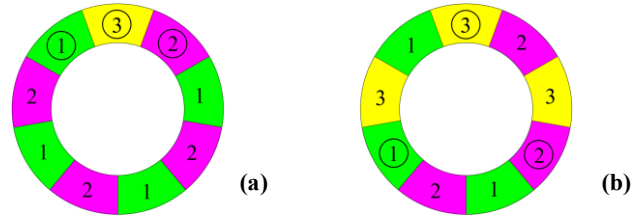

**Fig5. Two methods to color an odd cycle.**

There are two ways to color an odd cycle. One way is to color the cycle mainly with two colors except for one face for the third color, just as Figure 5(a). Another way is to partition the cycle into three odd paths with three faces and assign different colors to the three faces, then color each path with the two colors of its end faces. A hexagonal lattice around an odd face can be colored with similar ways.

**(3+1)-Coloring**

**Theorem 5.** The hexagonal lattice around an odd face can be colored mainly with three colors except for one congruence chain with a fourth color.

**Proof.** The lattice around an odd face can be unfolded along a congruence chain from the odd face as Figure 6(a). If the unfolded lattice is colored with 3-coloring, the adjacent faces along two sides of the unfolding chain are incompatible as they have the same color. These faces belong to two adjacent chains of the unfolding chain. If one of these two chains is colored with a fourth color as Figure 6(b), the conflict will be resolved. with this coloring, the selected chain and the unfolding chain form a Kempe chain separating the coloring on two sides.

□

In this paper, such coloring is called (3+1)-coloring. The unfolding chain is not necessarily linear, and so is the chain for the fourth color. Any congruence chain that goes outward from the neighborhood of the odd face can be selected for the fourth color.

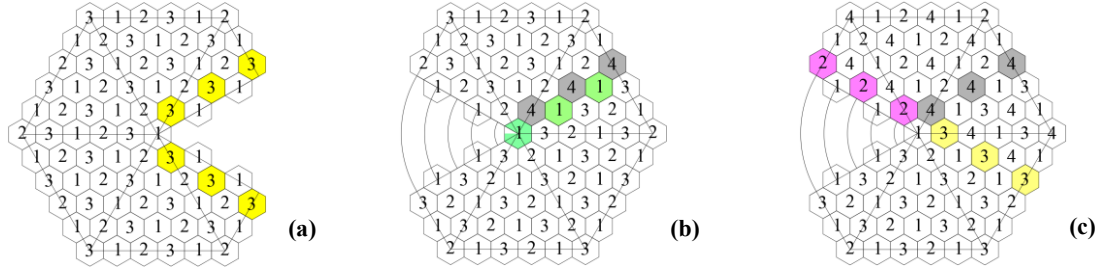

**Fig6. Face coloring around an odd face.**
(a). Incompatible coloring. (b). (3+1)-Coloring. (c). Multiple 3-coloring.

In the incompatible 3-coloring as Figure 6(a), the faces congruent to the odd face has no conflicts with 3-coloring. The conflicts are among the other faces. It is because on two sides of the unfolding line, the faces with two different colors are congruent across the line. With a walk around an odd face, these faces are no longer congruent to themselves. Therefore, those faces cannot be classified into two congruence classes. In this paper the congruence classes next to an odd face are all local congruence classes without a path around the face by default.

**Multiple 3-coloring**

**Theorem 6.** The hexagonal lattice around an odd face can be partitioned into three zones, each colored with a different 3-color subsets of 4 colors.

**Proof.** Just as mentioned above, the neighborhood of an odd face is an odd cycle, which can be partitioned into three odd paths with three faces for different colors, and the odd face itself can be colored with a fourth color. Selecting a congruence chain that goes outward from each partition face, the hexagonal lattice will be partitioned into three zones. As shown in Figure 6(c), each zone can be colored with 3 colors: the color of the odd face and colors of two partition faces. Each partition chain has an adjacent chain from the odd face and they together form a Kempe chain separating the coloring on two sides.

□

In this paper, such coloring is called multiple-3-coloring (M3-coloring for short). The (3+1)-coloring is the special case of the M3-coloring that some zones are shrunk.

The coloring for hexagonal lattice on a tube in Figure 3(b) is also (3+1)-coloring. As the path along the polyline is a cycle with length of the form $3n+1$, it can be partitioned into four paths whose length are all of the form $3n+1$. These paths can be colored with color cycles 1-2-3, 2-3-4, 3-4-1, and 4-1-2 sequentially. As the result, four partition faces will be colored with colors 1, 2, 3, and 4 separately. Now, if each linear congruence chain is colored with the color of corresponding face on the polyline, the lattice with be colored with M3-coloring.

## 5. Combinations of odd faces

### 5.1 Combination of a nice pair

Two odd faces can be combined as a nice pair if they are connected by a congruence chain. With sharing the partition chains, the hexagonal lattice around a nice pair can be colored with (3+1)-coloring as Figure 7(a), or M3-coloring as Figure 7(b). In these two colorings, the outskirt of the nice pair is colored with 3 colors.

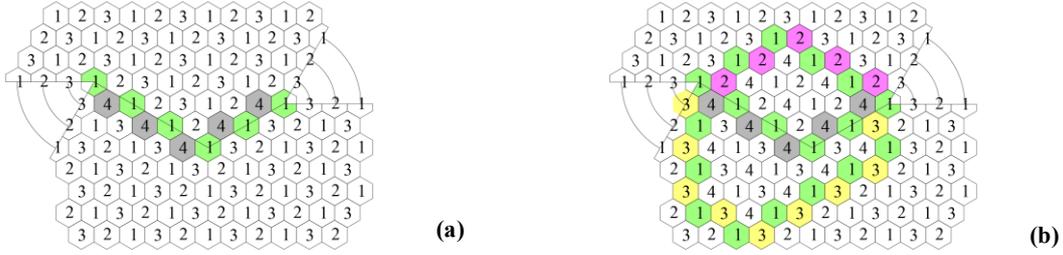

**Fig7. Combination of a nice pair for coloring.**

**Theorem 7.** The outskirt of two odd faces is 3-colorable if and only if they are a nice pair.

**Proof.** In the incompatible 3-coloring as Figure 6(a), the color of the congruence faces across the unfolding line can be described with cycle notation. In 3-coloring with colors 1, 2, and 3, the permutation around odd-face with color 1 is (2, 3). Similarly, the permutation around odd-face with color 2 is (1, 3), and the permutation around odd-face with color 3 is (1, 2). If the lattice around two odd faces is unfolded and colored as Figure 8, the color permutation across the common unfolding line is the product of permutations around two odd faces. Let face color of a nice pair is color 1 as in Figure 8(a), then the color change across the common unfolding line is $(2, 3) \cdot (2, 3) = (1)(2)(3)$. It means the face colors on the two sides are consistent. The outskirt of two odd faces is 3-colorable.

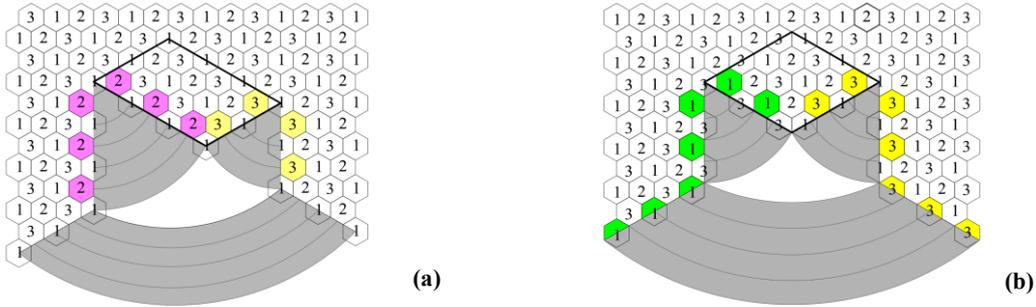

**Fig8. Coloring around two odd faces.**

Similarly, it is easy to get that if two odd faces are not a nice pair, the face colors across the common unfolding line are not consistent, and the outskirt is not 3-colorable, just as Figure 8(b).

□

As the outskirt is 3-colorable, a nice pair can be treated as an even face in further combination of other odd faces.

**Theorem 8.** The odd faces of a polyhedral graph can all be combined into nice pairs.

**Proof.** First, we can combine the adjacent odd faces suitable to nice pairs with priority for closer. If there are four odd faces of which any two adjacent ones are incongruent as in Figure 9(a), they can be combined as follows. As face *A* and *C* are not congruent to *B* and *D*, they are congruent to each other along a chain across *BD*. So, they can be combined into a nice pair. After this combination, the outskirt of face *A* and *C* can be colored with 3-coloring. Then face *B* and *D* have the same color and they can be combined as a nice pair, just as Figure 9(b).

The combined nice pairs can be treated as even faces and such combinations can go on. At last, there are at most two odd faces left since the total number is even. The outskirt of these two faces is also the outskirt of other nice pairs, which can be colored with 3-coloring. According to Theorem 7, they must be a nice pair.

□

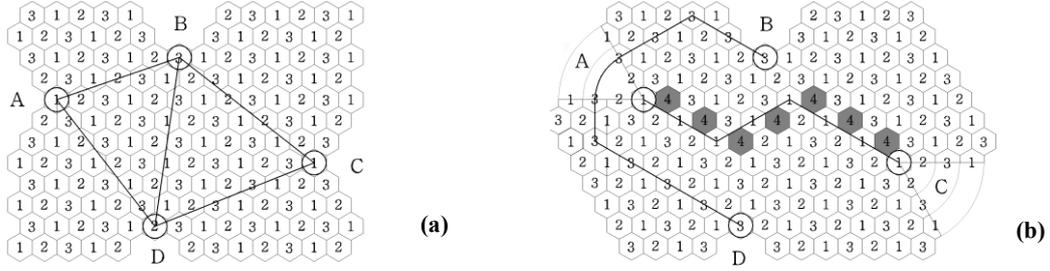

**Fig9. Combination of four incongruent odd faces into two nice pairs.**

### 5.2 Combination of an incongruent triplet

Three odd faces form an inconsistent triplet if they are incongruent with each other in their common neighborhood.

**Theorem 9.** An incongruent triplet can be combined as an equivalent odd face for M3-coloring, and the equivalent face is not congruent to these three faces.

**Proof.** Just as in Figure 10(a), odd faces *A, B,* and *C* are incongruent inside the triangle *ABC*. But along a curve around one of these faces, as shown in the figure, the other two faces are congruent to each other. Outside the triangle, there is a class of faces incongruent with *A, B,* and *C,* but congruent to face *A* if across *BC*, congruent to face *B* if across *AC*, and congruent to face *C* if across *AB*. These faces can be colored with a fourth color without conflict.

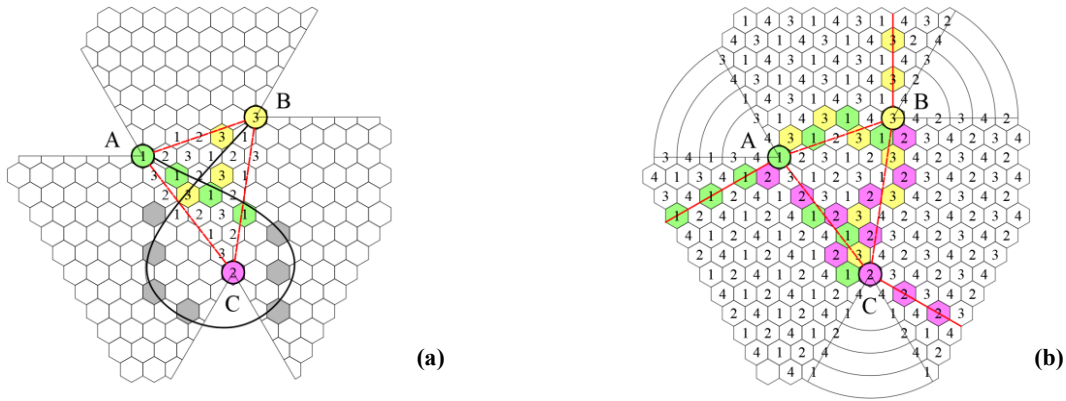

**Fig10. Combination of an incongruent triplet.**

Just as in Figure 10(b), with a congruence chain that goes outward from each odd face, the lattice outside the triangle will be partitioned into three zones, each can be colored with three colors: the colors of two odd faces on the corner and the fourth color. Such M3-coloring for the outskirt of an incongruent triplet is just like that around an odd face in Figure 6(c), and the equivalent face is not congruent to the three faces in the triplet.

□

### 5.3 Global Combination

A nice pair can be treated as an even face, and an incongruent triplet can be treated as an odd face for further combination. At last, according to Theorem 8 and 9, all odd faces can be combined into clusters whose outskirts can be colored with 3-coloring.

For example, four incongruent odd faces in Figure 9 also can be combined as follows. Faces *A*, *B*, and *D* can be combined into an incongruent triplet. As face *C* is incongruent with those three from the outside, face *C* and the equivalent face of the triplet are a nice pair. Therefore, the faces around these four odd faces can be colored as in Figure 11(a). with this coloring, there is a Kempe chain between each two of the four odd faces, partitioning the graph into four zones for M3-coloring.

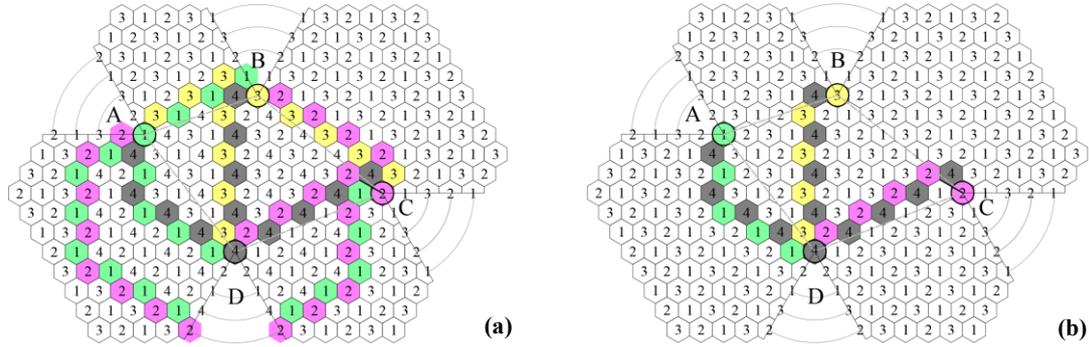

**Fig11. Combination of four incongruent odd faces.**

With three zones shrinking to the minimum, the M3-coloring in Figure 11(a) can be transformed to (3+1)-coloring in Figure 11(b), where only three congruence chains from face *D* toward the other three faces are colored with a fourth color.

As the outskirt can be colored with 3-coloring, the combination of such four incongruent faces can also be treated as an even face in further combination.

The combinations treated as even faces will be embedded in the zones of other combinations. A combination being embedded in a zone of another combination means they have a common zone, i.e. they are connected by a 3-colorable zone. The two clusters can be colored independently so long as their common zone is colored with 3-coloring. Therefore, with a global combination, a graph will be decomposed into independent subgraphs for M3-coloring.

# 6. Combination graph and graph triangulation

## 6.1 Combination graph

The combination and partition of graphs for M3-coloring can be illustrated with combination graphs in Figure 12. In the combination graphs, a circle marked with a number represents a node face and its color, and a line/curve segment connecting two circles represents a Kempe chain.

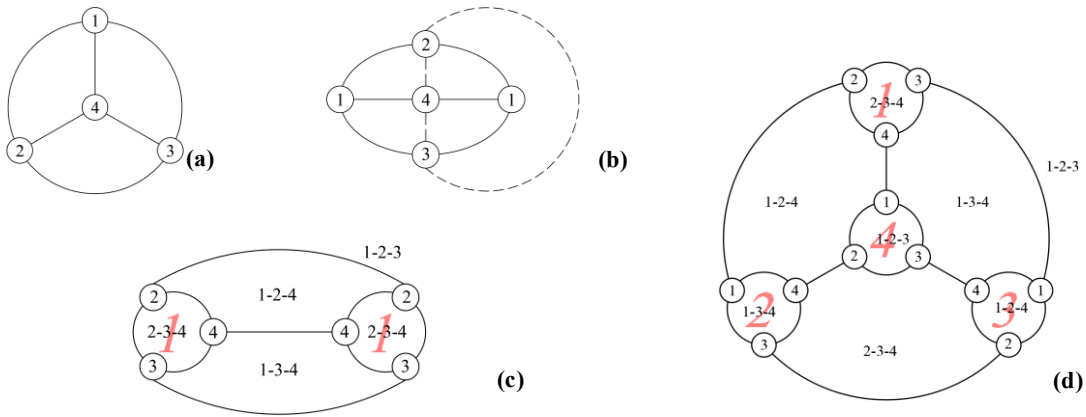

**Fig12. Examples of combination graphs.**

The combination of four incongruent odd faces in Figure 11(a) can be illustrated with the graph in Figure 12(a). In fact, this combination graph is the basic combination/partition for M3-coloring.

If three of the four faces are at infinity, it represents the partition around an odd face shown in Figure 6(c). If one of the four faces is at infinity, it represents the combination of an incongruent triplet shown in Figure 10(b).

The combination of a nice pair in Figure 7(b) can be illustrated with the graph in Figure 12(b). It can be seen as the splicing of two combinations of Figure 7(a). Three dashed line segments mean three Kempe chains that can be omitted as no use for coloring partition.

Figures 12(c), and 12(b) represent two combinations of equivalent odd faces. Figure 12(c) is for a nice pair and Figure 12(d) is for four incongruent faces. In these two combination graphs, the colors for equivalent faces (bold number) and colors for 3-coloring zones are annotated for clarity.

### 6.2    Combination adjustment

Not all the combinations of a graph can be successfully partitioned for M3-coloring. For example, many graphs cannot be partitioned for M3-coloring through all odd faces combined into nice pairs. After global combination, it is necessary to select reasonable Kempe chains for partition. To get reasonable Kempe chains, following adjustments can be used if necessary.

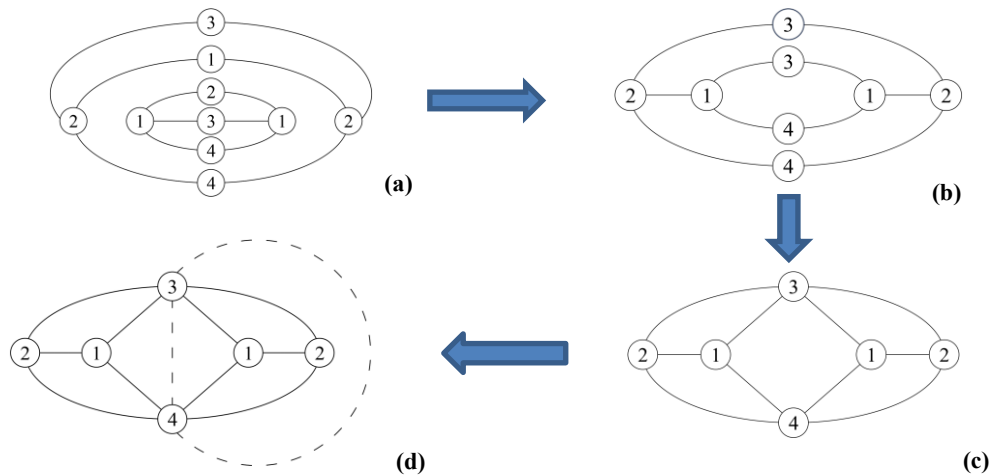

Fig13. Adjustment of Kempe chains.

**(1) Two adjacent Kempe chains with the same colors can be merged.**

Just as in Figure 13(a), a nice pair with color 1 is embedded in a zone of another nice pair with color 2. Each pair has a Kempe chain with colors 1 and 2 for partition. Such two Kempe chains can be merged. After merging, the common part can be omitted as the colorings on the two sides are the same. The remaining parts are two Kempe chains connecting the two pairs in Figure 13(b).

**(2) Two adjacent Kempe chains with only one common color can be joined at one face.**

Just as in Figure 13(b), two adjacent Kempe chains with common color 3 can have a common face of color 3, and so do two adjacent Kempe chains with common color 4. If so, we can get the result in Figure 13(c).

**(3) Partitioning a graph into two independent parts with Kempe chains.**

Just as the dashed line shown in Figure 13(d), two Kempe chains of colors 3 and 4 form a closed circle. the graph is partitioned into two parts which can be colored independently: the color of 1 and 2 can be exchanged in each part. If two parts are separated by Kempe chains between three nodes, just as dashed lines in Figure 12(b), the two parts can also be colored independently.

### 6.3  Graph triangulation

The classic partition of a graph for multiple 3-coloring is graph triangulation. Combination graphs in Figures 12(a) and 12(b) are certainly graph triangulations. Combination graphs in Figures 12(c) and 12(d) are also graph triangulations where equivalent odd faces are used as triangulation nodes. The result in Figure 13(d) is also a graph triangulation.

The "triangulation" here is in the topological sense. Generally, Kempe chains are not linear, and the zones are not triangular. As the three node faces of each zone are incongruent with each other, the face zone can be 3-colored with the node colors. With proper triangulation, if the combination graph can be 4-colored, the graph can be colored with M3-coloring.

In the graph triangulation, even faces can also be selected as node faces so long as they are incongruent with the adjacent nodes, just as in Figure 12(b). After triangulation, the number of partition Kempe chains connecting to an odd face node should be odd, and that to an even face node should be even.

In fact, each graph itself is a graph triangulation: any three adjacent faces are incongruent with each other, and they form a triangular zone. Proper triangulation can transform a graph into simple combination graphs for coloring, one for each independent subgraph.

To get proper triangulation, sometimes it is necessary to adjust the combinations. For example, two adjacent nice pairs may be merged and triangulated into a tetrahedron; a tetrahedron may be refined into a (3,6)-fullerene (a tetrahedron is a special kind of (3,6)-fullerene); the combinations of 12 odd faces in Figure 12(d) may be triangulated into an icosahedral fullerene.

Generally, independent subgraphs can be triangulated into symmetric structures, such as tetrahedrons (3,6)-fullerenes (including the tetrahedron), icosahedral fullerenes, etc., or their splicing as shown in figures 12(b) and 13(d). Those with asymmetric structures can be decomposed further or merged with other combinations. Certainly, the combination graphs may not be symmetric indeed because of the mix of odd faces and equivalent odd faces.

## 7.  Orbital coloring for symmetric graphs

As mentioned in section 2, all (3,6)-fullerenes are face R4-colorable. The lattice of a (3,6)-fullerene is just like the result of a three-axis weaving. The faces along the weaving lines of each axis form a set of orbifolds. As shown in Figure 14, with R4-coloring, each orbifold is colored with two colors, and its adjacent orbifolds are colored with another two. In this paper, such coloring is called orbital coloring. As each orbifold of a (3,6)-fullerene is either a path or an even cycle, there are no conflicts in the orbital coloring.

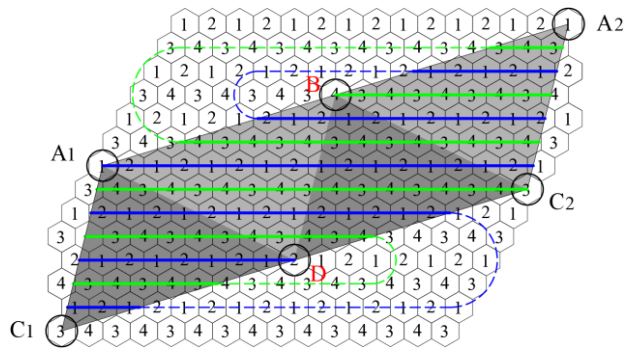

**Fig14. Orbit coloring of a (3,6)-fullerene.**

**Theorem 10.** Icosahedral fullerenes can be face-colored with orbital coloring.

**Proof.** An icosahedral fullerene has twenty identical equilateral triangular cells. The triangular cell has two parameters shown in Figure 15(a). The lattice lines from three vertices of a cell can converge in a point as point $D$ in the figure, so the cone manifolds with vertices around a cell can be unfolded through the lattice lines from one point.

If $q = 0$, the fullerene manifold can be unfolded as Figure 15(b). Two light-colored regions in the figure indicate the same region of the manifold. After unfolding, the faces along the horizontal line form a set of orbifolds. Because of the symmetry of the fullerene, these orbifolds are either a path or an even cycle. So, they can be colored with orbital coloring.

Similarly, if $p = q$, the fullerene can be unfolded as in Figure 15(c). If $p > q$ and $q \neq 0$, the fullerene can be unfolded as in Figure 15(d). It is easy to verify that in these two cases, the faces along the horizontal line in the unfolding also form orbifolds and can be colored with orbital coloring.

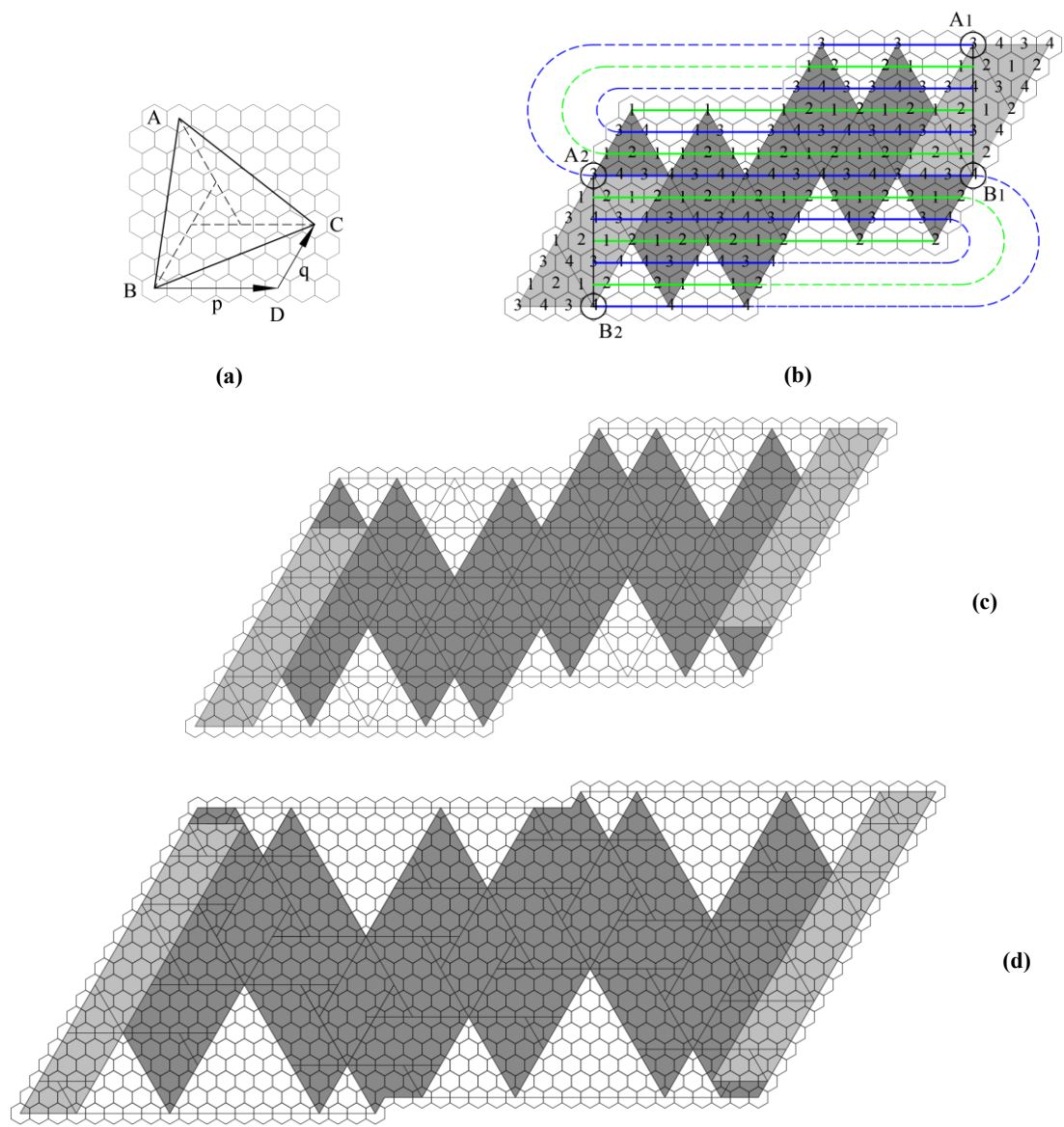

(a) (b) (c) (d)

**Fig15. Unfolding of icosahedral fullerenes for orbital coloring.**

□

# 8. Conclusion and outlook

Multiple 3-coloring is an ideal approach to the 4-coloring of planar graphs. Generally, a planar graph can be decomposed into independent subgraphs, each has a symmetric combination graph that can be colored with orbital coloring, and the subgraph can be colored with M3-coloring according to the orbital coloring of the combination graph. If so, the planar graph can be colored with 4 colors. If it can be proved that all planar graphs can be decomposed and colored with this method, the theoretical proof for the four-color theorem can be reached. To reach such a proof, the symmetries of planar graphs should also be considered. Just as fullerenes [10], each planar graph has certain intrinsic symmetries.